\documentclass{amsart}
\usepackage{amssymb}
\newtheorem{thm}{Theorem}[section]
\newtheorem{cor}[thm]{Corollary}
\newtheorem{lem}[thm]{Lemma}
\newtheorem{klem}[thm]{\bf {Key Lemma}}

\theoremstyle{definition}
\newtheorem{defn}[thm]{Definition}
\theoremstyle{remark}

\numberwithin{equation}{section}

\newcommand{\R}{\mathbb R}

\newcommand{\rt}{\rightarrow}

\newcommand{\C}{{\mathbb C}^n}
\newcommand{\B}{{\mathbb B}^n}

\begin{document}

\title[Negative curvature on products]{Negative sectional curvature and the \\ product complex structure }
\author{Harish Seshadri}
\address{Department of Mathematics,
Indian Institute of Science, Bangalore 560012, India}
\email{harish@math.iisc.ernet.in}

 \thanks{This work was done under DST grant SR/S4/MS:307/05}

\begin{abstract}
Let $M=M_1 \times M_2$ be a product of complex manifolds. We
prove that $M$ cannot admit a complete K\"ahler metric with
sectional curvature $K<c<0$ and Ricci curvature $Ric
> d$, where $c$ and $d$ are constants.

In particular, a product domain in $\C$ cannot cover a compact
K\"ahler manifold with negative sectional curvature.

On the other hand, we observe that there are complete K\"ahler
metrics with negative sectional curvature on $\C$. Hence the
upper sectional curvature bound is necessary.

\end{abstract}
\maketitle
\section{Introduction}
The interplay between the curvature and the underlying complex
structure of a K\"ahler manifold is a central theme in complex
differential geometry. In this article we prove a general result
ruling out the existence of negatively curved K\"ahler metrics on
product complex manifolds. The inspiration for our result is the
classical Preissmann theorem stating that the fundamental group
of a compact negatively curved Riemannian manifold does not
contain ${\mathbb Z} \oplus {\mathbb Z}$ as a subgroup. In fact,
if the factors are assumed compact, then our theorem follows from
Preissmann's theorem. In the general situation the obstruction to
negative curvature arises from the complex structure rather than
the topology of $M$.
\begin{thm}\label{mai}
Let $M=M_1 \times M_2$ be a product of complex manifolds $M_1$
and $M_2$ with $dim \ M_i \ge 1$, $i=1,2$. Then $M$ cannot admit a
complete K\"ahler metric with sectional curvature $K<c<0$ and
Ricci curvature $Ric > d$, where $c$ and $d$ are constants.
\end{thm}

An important feature of Theorem \ref{mai} is that no assumptions
are made about the factors $M_i$. Let us compare our result with
the work of P. Yang ~\cite{yan} and F. Zheng ~\cite{zhe} which
again consider the interaction between the product complex
structure and negative curvature. On the one hand, both these
papers assume only negative or nonpositive {\it holomorphic
bisectional} curvature $B$. On the other hand, the assumptions on
the factors $M_i$ are more stringent. Yang's paper rules out the
existence of complete K\"ahler metrics with $d<B<c<0$ on {\it
polydiscs} (more generally, bounded symmetric domains of rank
$>1$) while Zheng classifies all metrics with $B \le 0$ on
products of {\it compact} manifolds. \vspace{2mm}

Let us consider Theorem \ref{mai} in the context of a basic
question regarding negatively curved K\"ahler manifolds: Is every
simply-connected, complete K\"ahler manifold $M$ with sectional
curvatures bounded between two negative constants biholomorphic
to a bounded domain in $\C$ ? (cf. ~\cite{ao}, ~\cite{wu} and
~\cite{gre}). This question is still open, even in the special
case of $M$ being the universal cover of a compact K\"ahler
manifold with negative sectional curvature. In this case if one
imposes further restrictions on $M$, there are interesting
results due to B. Wong ~\cite{won}, J.-P. Rosay ~\cite{ros} and
S. Frankel ~\cite{fra}. In these works $M$ is only assumed to be
the universal cover of a compact complex manifold. The Wong-Rosay
theorem implies that if such an $M$ is a domain in $\C$ with
$C^2$-smooth boundary, then $M$ has to be biholomorphic to $\B$.
According to the work of Frankel, if $M$ is a bounded convex
domain, then $M$ has to be biholomorphic to a bounded symmetric
domain. As a corollary of our theorem, we have \vspace{2mm}

\begin{cor}
A product domain in $\C$ cannot cover a compact K\"ahler manifold
with negative sectional curvature. In fact, product domains do
not admit complete K\"ahler metrics with pinched negative
sectional curvature.
\end{cor}

{\bf Remark:} Theorem \ref{mai} has other intriguing
complex-analytic implications. For instance, it follows that the
unit ball in $\C$ is not biholomorphic to a product of complex
manifolds, a fact which is probably known. \vspace{2mm}

Regarding the necessity of the assumptions on the curvature in
Theorem \ref{mai}, we show in the last section that the
calculations in ~\cite{kle} can be adapted to get a complete
K\"ahler metric on $\C$ with sectional curvature $d<K<0$. Hence
the upper bound $K<c<0$ is necessary in Theorem \ref{mai}.
However, the following question is still open: Is Theorem
\ref{mai} still valid if we drop the lower bound on $Ric$ ? Also,
it would be interesting to know if the result remains true if one
replaces sectional curvature with holomorphic bisectional
curvature. \vspace{2mm}

The proof of Theorem \ref{mai} is ``soft", modulo the use of Yau's
Schwarz Lemma. By this we mean that only ``coarse" geometric
ideas in the sense of Gromov are used. It can be summarized as
follows: Without loss of generality, one can assume that $M$ is
simply-connected. The key observation is that by using Yau's
Schwarz Lemma one can show that $(M,g)$ has to be bi-Lipschitz to
a product of non-compact Riemannian manifolds $(M_1, g_b) \times
(M_2,g_a)$. Finally we prove that such a product cannot be
Gromov-hyperbolic. But $(M,g)$ being a simply-connected
Riemannian manifold of negative sectional curvature (bounded away
from zero) has to be Gromov-hyperbolic. \vspace{2mm}

{\bf Acknowledgement:} I would like to thank Fangyang Zheng for his helpful comments 
regarding this paper.

\section{proof}
If $(X,d), \ (X_1,d_1)$ and $(X_2,d_2)$ are metric spaces, we
write $(X,d)=(X_1,d_1) \times (X_2,d_2)$ to mean that $X=X_1
\times X_2$ and $d(p,q)=d_1(p_1,q_1)+d_2(p_2,q_2)$ for all
$p=(p_1,p_2), \ q=(q_1,q_2)$ in $X$.

 For $0<L< \infty$, metrics \ $d_1$ and $d_2$ on $X$ are said to be {\it L-bi-Lipschitz} (or just {\it bi-Lipschitz}) if
$$L^{-1} \ d_1(x,y) \le d_2(x,y) \le L \ d_1(x,y)$$
for all $x,y \in X$.

Let $M=M_1 \times M_2$ be a complex manifold with K\"ahler metric
$g$. For $a \in M_1$, $b \in M_2$, let $i_b: M_1 \rt M_1 \times
M_2$ and $i_a:M_2 \rt M_1 \times M_2$ be $i_b(x)=(x,b)$ and
$i_a(y)=(a,y)$. Also, let $g_a:= i^\ast_a (g)$ and
$g_b:=i^\ast_b(g)$ denote the induced K\"ahler metrics on $M_2$
and $M_1$ respectively.
\begin{klem}\label{bi}
Let $M=M_1 \times M_2$ be a product complex manifold. Suppose that
$g$ is a complete K\"ahler metric on $M$ with $Ric>d$ and
$H<c<0$, where $H$ denotes holomorphic sectional curvature.

Then $(M,d)$ is bi-Lipschitz to $(M_1,d_b) \times (M_2, d_a)$,
where $d, \ d_b$ and $d_a$ are the distance functions associated
to $g,\ g_b$ and $g_a$ respectively. Here $a \in M_1$ and $b \in
M_2$ are arbitrary.

\end{klem}
\begin{proof}
We recall Yau's Schwarz Lemma which is the main ingredient in the
proof :

\begin{thm}\label{yau}{\rm (Yau's Schwarz Lemma ~\cite{yau})}
Let $(M,g)$ be a complete K\"ahler manifold and with Ricci
curvature $Ric \ge d$ and let $(N,h)$ a Hermitian manifold with
holomorphic sectional curvature $H \le c<0$, where $c$ and $d$ are
constants. If $f: M \rt N$ is a non-constant holomorphic map, then
$d <0$ and $f^\ast (h) \le \frac {d}{c} g$.

Hence $d_h(f(p),f(q)) \le {\sqrt \frac {d}{c}} \ d_g(p,q)$ for any
$p,q \in M$.

\end{thm}

We note that one only needs an upper bound on the {\it holomorphic
sectional} curvature $H$, as opposed to the full sectional
curvature, of $N$ to apply this theorem (in Yau's original paper
~\cite{yau}, an upper bound on holomorphic bisectional curvature
of $N$ is assumed. For a proof involving only $H$, see
~\cite{zhe2}. ).

Now we restrict to $M=M_1 \times M_2$. Observe that the
holomorphic sectional curvatures of $g_a$ and $g_b$ are less than
$c$, \ since ${\{a\} \times M_2}$ and ${M_1 \times \{b\}}$ are
complex submanifolds of $M$.

Let $\pi_i:M \rt M_i$ denote the projections. Let $x,y \in M_1$
and $b,b' \in M_2$. Applying Theorem \ref{yau} to $ \pi_1:(M,g)
\rt (M_1,g_b),$ we conclude that
$$  d_b( x, y) \le L \  d((x,b'),(y,b')),$$
where $L = \sqrt {\frac {d}{c}}$. Since $g_{b'}$ is induced from
$g$, we have $d((x,b'),(y,b')) \le d_{b'}(x,y)$ and hence
$$ d_b(x,y) \le L  \ d_{b'}(x,y).$$
Interchanging $b$ and $b'$, we see that the distance functions
$d_b$ and $d_{b'}$ are L-bi-Lipschitz on $M_1$ for any $b, b' \in
M_1$. Similarly $d_a$ and $d_{a'}$ are L-bi-Lipschitz on $M_2$ for
any $a, a' \in M_1$. \vspace{2mm}

Fix $a \in M_1$ and $b \in M_2$. Let $p=(p_1,p_2), \ q=(q_1,q_2)$
be arbitrary points in $M$. Again, by Theorem \ref{yau} applied to
$\pi_1: (M,g) \rt (M_1, g_{q_2})$, we get
$$ d_{q_2}(p_1,q_1) \le L \ d (p,q).$$
Similarly $d_{p_1}(p_2,q_2) \le L \ d (p,q)$. Adding these two
inequalities and using the L-bi-Lipschitz equivalence of all
metrics on $M_1$ and $M_2$ (as proved in the previous paragraph),
we get
\begin{equation}\label{one}
d_{b}(p_1,q_1) + d_{a}(p_2,q_2) \le 2L^2 \ d(p,q).
\end{equation}
On the other hand, letting $r=(p_1,q_2)$ and applying the
triangle inequality we have
\begin{align} \label{two}
d(p,q) & \le d(p,r)+d(r,q) \\ \notag
        & \le d_{p_1}(p_2,q_2) +d_{q_2}(p_1,q_1) \\ \notag
        & \le L \ (d_a(p_2,q_2)+d_b(p_1,q_1)).
\end{align}
In the second inequality we have again used the fact that if $p,q
\in M_1 \times \{b\} $ for some $b \in M_2$, then $d_b(p,q) \ge
d(p,q)$, with similar inequality holding for $p,q \in \{a\} \times
M_2$.

Combining (\ref{one}) and (\ref{two}), we see that $(M,d)$ is
bi-Lipschitz to $(M_1, d_a) \times (M_2,d_b)$.

\end{proof}

{\bf Remark}: The use of Yau's Schwarz Lemma in this proof was
inspired by Yang's article ~\cite{yan}. \vspace{4mm}

 Let us now recall the concept of {\it
Gromov-hyperbolicity} of metric spaces. For the sake of clarity,
we give definitions that are more general than are actually
needed. For further details, we refer the reader to ~\cite{ghy}.

\begin{defn}
Let $(X,d)$ be a metric space. A {\it geodesic} between $ p,q \in
X$ is an isometric map $\gamma: [0, d(p,q)] \rt X$ with $\gamma
(0)=p, \ \gamma (d(p,q))=q$. A metric space is {\it geodesic} if
any there is a geodesic between any two points in $X$.

A {\it geodesic triangle} $\Delta$ in $X$ is a union of images of
three geodesics $\gamma_i: [a_i,b_i] \rt X $ with $ \gamma_i(b_i)=
\gamma_{i+1}(a_{i+1})$, where $i$ is taken mod 3. The image of
each $\gamma_i$ is called a {\it side} of $\Delta$.

A geodesic metric space $X$ is {\it Gromov-hyperbolic} or {\it
$\delta$-hyperbolic} if there exists $\delta >0$ such that for
any geodesic triangle $\Delta$ in $X$, the $\delta$-neighbourhood
of any two sides in $\Delta$ contains the third side.

\end{defn}
An important class of Gromov-hyperbolic spaces are the
simply-connected, complete Riemannian manifolds with sectional
curvature bounded above by a negative constant. The
Gromov-hyperbolicity of these is a consequence of Toponogov's
comparison theorem.




The following basic lemma captures the coarse geometric nature of
Gromov-hyperbolicity. This lemma is actually valid under the
weaker hypothesis of {\it quasi-isometric} equivalence, but we
only require the bi-Lipschitz case.
\begin{lem}\label{also}{\rm (cf. ~\cite{ghy})}
Let $d$ and $d'$ be bi-Lipschitz metrics on $X$. If $(X,d)$ is
Gromov-hyperbolic, then so is $(X,d')$.
\end{lem}
For the next lemma, recall that a {\it geodesic ray} in $X$ is an
isometric map $\gamma:[0, \infty) \rt X$.
\begin{lem}\label{not}
Let $(X,d)=(X_1,d_1) \times (X_2,d_2)$. If there is a geodesic
ray in each $(X_i,d_i)$ for i=1,2, then $(X,d)$ is not
Gromov-hyperbolic.
\end{lem}
\begin{proof}
For $i=1,2$, let $\gamma_i$ be a geodesic ray in $(X_i,d_i)$. For
each $n \in {\mathbb Z}^+$, define the map $\sigma_n: [0,2n] \rt
X$ by
\begin{align} \notag
\sigma_n (t) & = ( \gamma_1(t), \gamma_2 (n)) \ \ \ \ \ \ \ \ \
{\rm for} \ \ 0 \le t \le n  \\ \notag
          & = (\gamma_1(n), \gamma_2(2n-t)) \ \ \ \ {\rm for} \ \ n \le t \le 2n
\notag
\end{align}
Then $\sigma_n: [0,2n] \rt (X,d)$ is a geodesic and if $S_1
=\gamma_1 ([0,n]) \times \{\gamma_2(0)\}$, $S_2 =\{\gamma_1 (0)
\} \times \gamma_2([0,n])$ then $\Delta_n= S_1 \cup S_2 \cup
\sigma_n([0,2n])$ is a geodesic triangle. But the distance between
$\sigma_n(n)=(\gamma_1 (n), \gamma_2(n))$ and $S_1 \cup S_2$ is
at least $n$. Hence $(X,d)$ cannot be $\delta$-hyperbolic for any
$\delta$.

\end{proof}

Now we can complete the proof of Theorem \ref{mai}: By taking the
universal cover of $M$ if necessary, we can assume that $M$ is
simply-connected. As noted earlier, $M$ will then be
Gromov-hyperbolic.
 On the other hand, since $M$ is diffeomorphic to $\R^{2n}$ by
 the theorem of Cartan-Hadamard, both $M_1$ and $M_2$ are
 noncompact. Moreover, since $(M,g)$ is complete, $(M_1,g_b)$ and
 $(M_2,g_a)$ are complete for any $a \in M_1$ and $b \in M_2$. Now
 one can always find a geodesic ray in any complete noncompact Riemannian
 manifold. By Lemma \ref{not}, $(M_1,d_b) \times (M_2,d_a)$ is
 {\it not} Gromov-hyperbolic. Hence, neither is $(M,d)$ by Lemma
 \ref{bi} and Lemma \ref{also}. This contradiction completes the proof.  \hfill Q.E.D.

 \section{Negatively curved K\"ahler metrics on $\C$}
 The example in this section is adapted from ~\cite{kle} and we
 refer the reader to it for further details. It is clear from
 Klembeck's computations that if
 $$g =  \sum_{i,j=1}^n \frac {\partial f(r^2) }{ \partial z_i \partial \bar z_j}
 \ dz_i \otimes d \bar z_j, $$
where $f: \R \rt \R$ and $r^2 =\vert z_1 \vert^2+...+ \vert z_n
\vert^2$, then the following are sufficient conditions for $g$ to
be a complete K\"ahler metric of strictly negative sectional
curvature on $\C$:

(a) \ $f'(r^2)+r^2f''(r^2)>0$,

(b) \ $\int_o^\infty \sqrt{f'(r^2)+r^2f''(r^2)} dr = \infty$,

(c) \ $f''(r^2)>0$,

(d) \ $f''(r^2)+r^2f'''(r^2)-r^2 \frac {f''(r^2)^2}{f'(r^2)} > 0$,

(e) $ \frac {1}{r} \  \frac {\partial}{\partial r} \Big( r
\frac{\partial}{\partial r} ln \Big(f'(r^2)+r^2f''(r^2) \Big)
\Big) >0$

If we let $f(x)=e^x$, then it can be checked that all the above
conditions are satisfied. Moreover, the sectional curvatures of
$g$ will be bounded below. Indeed, some sectional curvatures will
tend to zero exponentially fast while some will equal $-2$ as $r
\rt \infty$.

Hence the upper sectional curvature bound in Theorem \ref{mai} is
necessary.

\end{document}